\documentclass[12pt]{amsart}
\usepackage{amssymb,latexsym}
\usepackage[usenames,dvipsnames]{pstricks}
\usepackage{epsfig}
\usepackage{pst-grad} 
\usepackage{multicol}

\begin{document}
\title[Exterior algebra]{The exterior algebra and central notions in
mathematics}
\author{Gunnar Fl{\o}ystad}
\address{Matematisk Institutt\\
         Johs. Brunsgt. 12\\ 
        5008 Bergen} 
\email{gunnar@mi.uib.no}

\date{\today}

\maketitle

\begin{center}
{\it Dedicated to Stein Arild Str{\o}mme (1951-2014)}
\end{center}

\medskip
\begin{quote}
``The neglect of the exterior algebra is the mathematical tragedy
of our century.''
\vskip 1.5mm
-Gian-Carlo Rota, \emph{Indiscrete thoughts (1996)}
\end{quote} 


\theoremstyle{plain}
\newtheorem{theorem}{Theorem}[section]
\newtheorem{corollary}[theorem]{Corollary}
\newtheorem*{main}{Main Theorem}
\newtheorem{lemma}[theorem]{Lemma}
\newtheorem{proposition}[theorem]{Proposition}

\theoremstyle{definition}
\newtheorem{definition}[theorem]{Definition}

\theoremstyle{remark}
\newtheorem{notation}[theorem]{Notation}
\newtheorem{remark}[theorem]{Remark}
\newtheorem{example}[theorem]{Example}
\newtheorem{claim}{Claim}


\newcommand{\psp}[1]{{{\bf P}^{#1}}}
\newcommand{\psr}[1]{{\bf P}(#1)}
\newcommand{\op}{{\mathcal O}}
\newcommand{\opw}{\op_{\psr{W}}}
\newcommand{\go}{\op}

\newcommand{\ini}[1]{\text{in}(#1)}
\newcommand{\gin}[1]{\text{gin}(#1)}
\newcommand{\kr}{{\Bbbk}}
\newcommand{\pd}{\partial}
\newcommand{\vardel}{\partial}
\renewcommand{\tt}{{\bf t}}


\newcommand{\coh}{{{\text{{\rm coh}}}}}


\newcommand{\modv}[1]{{#1}\text{-{mod}}}
\newcommand{\modstab}[1]{{#1}-\underline{\text{mod}}}

\newcommand{\sut}{{}^{\tau}}
\newcommand{\sumit}{{}^{-\tau}}
\newcommand{\til}{\thicksim}

\newcommand{\totp}{\text{Tot}^{\prod}}
\newcommand{\dsum}{\bigoplus}
\newcommand{\dprod}{\prod}
\newcommand{\lsum}{\oplus}
\newcommand{\lprod}{\Pi}

\newcommand{\La}{{\Lambda}}

\newcommand{\sirstj}{\circledast}

\newcommand{\she}{\EuScript{S}\text{h}}
\newcommand{\cm}{\EuScript{CM}}
\newcommand{\cmd}{\EuScript{CM}^\dagger}
\newcommand{\cmri}{\EuScript{CM}^\circ}
\newcommand{\cler}{\EuScript{CL}}
\newcommand{\clerd}{\EuScript{CL}^\dagger}
\newcommand{\clerri}{\EuScript{CL}^\circ}
\newcommand{\gor}{\EuScript{G}}
\newcommand{\gF}{\mathcal{F}}
\newcommand{\gG}{\mathcal{G}}
\newcommand{\gM}{\mathcal{M}}
\newcommand{\gE}{\mathcal{E}}
\newcommand{\gD}{\mathcal{D}}
\newcommand{\gI}{\mathcal{I}}
\newcommand{\gP}{\mathcal{P}}
\newcommand{\gK}{\mathcal{K}}
\newcommand{\gL}{\mathcal{L}}
\newcommand{\gS}{\mathcal{S}}
\newcommand{\gC}{\mathcal{C}}
\newcommand{\gO}{\mathcal{O}}
\newcommand{\gJ}{\mathcal{J}}
\newcommand{\gU}{\mathcal{U}}
\newcommand{\mm}{\mathfrak{m}}

\newcommand{\dlim} {\varinjlim}
\newcommand{\ilim} {\varprojlim}

\newcommand{\CM}{\text{CM}}
\newcommand{\Mon}{\text{Mon}}


\newcommand{\Kom}{\text{Kom}}


\newcommand{\EH}{{\mathbf H}}
\newcommand{\res}{\text{res}}
\newcommand{\Hom}{\text{Hom}}
\newcommand{\inhom}{{\underline{\text{Hom}}}}
\newcommand{\Ext}{\text{Ext}}
\newcommand{\Tor}{\text{Tor}}
\newcommand{\ghom}{\mathcal{H}om}
\newcommand{\gext}{\mathcal{E}xt}
\newcommand{\id}{\text{{id}}}
\newcommand{\im}{\text{im}\,}
\newcommand{\codim} {\text{codim}\,}
\newcommand{\resol}{\text{resol}\,}
\newcommand{\rank}{\text{rank}\,}
\newcommand{\lpd}{\text{lpd}\,}
\newcommand{\coker}{\text{coker}\,}
\newcommand{\supp}{\text{supp}\,}
\newcommand{\Ad}{A_\cdot}
\newcommand{\Bd}{B_\cdot}
\newcommand{\Fd}{F_\cdot}
\newcommand{\Gd}{G_\cdot}


\newcommand{\sus}{\subseteq}
\newcommand{\sups}{\supseteq}
\newcommand{\pil}{\rightarrow}
\newcommand{\vpil}{\leftarrow}
\newcommand{\rpil}{\leftarrow}
\newcommand{\lpil}{\longrightarrow}
\newcommand{\inpil}{\hookrightarrow}
\newcommand{\pils}{\twoheadrightarrow}
\newcommand{\projpil}{\dashrightarrow}
\newcommand{\dotpil}{\dashrightarrow}
\newcommand{\adj}[2]{\overset{#1}{\underset{#2}{\rightleftarrows}}}
\newcommand{\mto}[1]{\stackrel{#1}\longrightarrow}
\newcommand{\vmto}[1]{\stackrel{#1}\longleftarrow}
\newcommand{\mtoelm}[1]{\stackrel{#1}\mapsto}

\newcommand{\eqv}{\Leftrightarrow}
\newcommand{\impl}{\Rightarrow}

\newcommand{\iso}{\cong}
\newcommand{\te}{\otimes}
\newcommand{\into}[1]{\hookrightarrow{#1}}
\newcommand{\ekv}{\Leftrightarrow}
\newcommand{\equi}{\simeq}
\newcommand{\isopil}{\overset{\cong}{\lpil}}
\newcommand{\equipil}{\overset{\equi}{\lpil}}
\newcommand{\ispil}{\isopil}
\newcommand{\vvi}{\langle}
\newcommand{\hvi}{\rangle}
\newcommand{\susneq}{\subsetneq}
\newcommand{\sgn}{\text{sign}}


\newcommand{\xd}{\check{x}}
\newcommand{\ortog}{\bot}
\newcommand{\tL}{\tilde{L}}
\newcommand{\tM}{\tilde{M}}
\newcommand{\tH}{\tilde{H}}
\newcommand{\tvH}{\widetilde{H}}
\newcommand{\tvh}{\widetilde{h}}
\newcommand{\tV}{\tilde{V}}
\newcommand{\tS}{\tilde{S}}
\newcommand{\tT}{\tilde{T}}
\newcommand{\tR}{\tilde{R}}
\newcommand{\tf}{\tilde{f}}
\newcommand{\ts}{\tilde{s}}
\newcommand{\tp}{\tilde{p}}
\newcommand{\tr}{\tilde{r}}
\newcommand{\tfst}{\tilde{f}_*}
\newcommand{\empt}{\emptyset}
\newcommand{\bfa}{{\bf a}}
\newcommand{\bfb}{{\bf b}}
\newcommand{\bfd}{{\bf d}}
\newcommand{\bfl}{{\bf \ell}}
\newcommand{\la}{\lambda}
\newcommand{\bfen}{{\mathbf 1}}
\newcommand{\ep}{\epsilon}
\newcommand{\en}{r}
\newcommand{\tu}{s}

\newcommand{\ome}{\omega_E}

\newcommand{\bevis}{{\bf Proof. }}
\newcommand{\demofin}{\qed \vskip 3.5mm}
\newcommand{\nyp}[1]{\noindent {\bf (#1)}}
\newcommand{\demo}{{\it Proof. }}
\newcommand{\demodone}{\demofin}
\newcommand{\parg}{{\vskip 2mm \addtocounter{theorem}{1}  
                   \noindent {\bf \thetheorem .} \hskip 1.5mm }}

\newcommand{\lcm}{{\text{lcm}}}


\newcommand{\dl}{\Delta}
\newcommand{\cdel}{{C\Delta}}
\newcommand{\cdelp}{{C\Delta^{\prime}}}
\newcommand{\dlst}{\Delta^*}
\newcommand{\Sdl}{{\mathcal S}_{\dl}}
\newcommand{\lk}{\text{lk}}
\newcommand{\lkd}{\lk_\Delta}
\newcommand{\lkp}[2]{\lk_{#1} {#2}}
\newcommand{\del}{\Delta}
\newcommand{\delr}{\Delta_{-R}}
\newcommand{\dd}{{\dim \del}}

\renewcommand{\aa}{{\bf a}}
\newcommand{\bb}{{\bf b}}
\newcommand{\cc}{{\bf c}}
\newcommand{\xx}{{\bf x}}
\newcommand{\yy}{{\bf y}}
\newcommand{\zz}{{\bf z}}
\newcommand{\mv}{{\xx^{\aa_v}}}
\newcommand{\mF}{{\xx^{\aa_F}}}

\newcommand{\Symm}{\text{Sym}}
\newcommand{\pnm}{{\bf P}^{n-1}}
\newcommand{\opnm}{{\go_{\pnm}}}
\newcommand{\ompnm}{\omega_{\pnm}}

\newcommand{\pn}{{\bf P}^n}
\newcommand{\hele}{{\mathbb Z}}
\newcommand{\nat}{{\mathbb N}}
\newcommand{\rasj}{{\mathbb Q}}

\newcommand{\dt}{\bullet}
\newcommand{\st}{\hskip 0.5mm {}^{\rule{0.4pt}{1.5mm}}}              
\newcommand{\disk}{\scriptscriptstyle{\bullet}}

\newcommand{\cF}{F_\dt}
\newcommand{\pol}{f}

\newcommand{\Rn}{{\mathbb R}^n}
\newcommand{\An}{{\mathbb A}^n}
\newcommand{\frg}{\mathfrak{g}}
\newcommand{\PW}{{\mathbb P}(W)}

\def\CC{{\mathbb C}}
\def\GG{{\mathbb G}}
\def\ZZ{{\mathbb Z}}
\def\NN{{\mathbb N}}
\def\RR{{\mathbb R}}
\def\OO{{\mathbb O}}
\def\QQ{{\mathbb Q}}
\def\VV{{\mathbb V}}
\def\PP{{\mathbb P}}
\def\EE{{\mathbb E}}
\def\FF{{\mathbb F}}
\def\AA{{\mathbb A}}

\section*{Introduction}
This note surveys how the exterior algebra and deformations or quotients
of it, gives rise to centrally important notions in five domains of 
mathematics: 
\begin{multicols}{2}
\begin{itemize}
\item Combinatorics
\item Topology
\item Lie theory
\item Mathematical physics
\item Algebraic geometry
\end{itemize}
\end{multicols}

The exterior algebra originated in the work of Herman
Grassman (1809-1877) in his book "Ausdehnungslehre" 
from 1844, and the thoroughly revised 1862 version, which
now exists in an English translation, \cite{Gra} from 2000.
Grassmann worked as professor at the gymnasium in 
Stettin, then Germany. Partly because being an original thinker,
and maybe partly because his education had not focused much on mathematics,
the first edition of the book had a more philosophical than mathematical
form, and therefore gained little influence in the mathematical community.
The second 1862 version was strictly mathematical. Nevertheless
it also gained little influence, perhaps because it had swung too far
to the other side and was scarce of motivation. Over 400 pages it developed the
exterior and interior product, and the somewhat lesser known regressive
product on the exterior algebra, which intuitively corresponds
to intersection of linear spaces. It relates this to geometry
and it also shows how analysis may be extended to functions of extensive
quantities. Only in the last two decades of the 1800's did
publications inspired by Grassmann's work achieve a certain mass.
It may have been with some regret that Grassmann
in his second version had an exclusively mathematical form, since he
in the foreword says "[extension theory] is not simply one among the other 
branches
of mathematics, such as algebra, combination theory or function theory,
but rather surpasses them, in that all fundamental elements are unified
under this branch, which thus as it were forms the keystone of the
entire structure of mathematics."

The present note indicates that he was not quite off the mark here. We do not
make any further connections to Grassmann's original presentation, but
rather present the exterior algebra in an entirely modern setting.
For more on the historical context of Grassmann, 
see the excellent history of vector analysis \cite{Crowe}, as well
as proceedings from conferences on Grassmann's manyfaceted legacy 
\cite{Schub} and \cite{Petsche}.
The last fifteen years have also seen a flurry of books advocating  
the very effective 
use of the exterior algebra and its derivation, the Clifford algebra, in 
physics, engineering and computer science. In the last section we report
briefly on this.


\section{The exterior algebra}

\subsection{Concrete definition}
Given a set $\{e_1, \ldots, e_n \}$ with $n$ elements. Consider
the $2^ n$ expressions $e_{i_1} \wedge e_{i_2} \wedge \ldots \wedge e_{i_r}$
(here $\wedge$ is just a place separator), where the $i_1, i_2, \ldots, i_r$
are strictly increasing subsequences of $1,2,\ldots, n$. 
From this we form the vector space $E(n)$ over a field $\kr$
with these expressions as basis elements. 

\begin{example}
When $n = 3$, the following eight expressions form a basis for $E(3)$:
\[ 1,  \, \,\, e_1, e_2, e_3, \,\, \, e_1 \wedge e_2, e_1 \wedge e_3, 
e_2 \wedge e_3, \,\,\, e_1 \wedge e_2 \wedge e_3. \]
\end{example}

The element $e_{i_1} \wedge e_{i_2} \wedge \ldots \wedge e_{i_r}$
is considered to have grade $r$, so we get a graded vector space $E(n)$.
Now we equip this vector space with a multiplication which we also
denote by $\wedge$. The basic rules for this multiplication are
\[ \mbox{i) }\quad  e_i \wedge e_i = 0, \qquad \mbox{ii) } \quad e_i \wedge 
e_j = - e_j \wedge e_i. \]
These rules together with the requirement that $\wedge$ is associative, 
i.e. 
\[ \mbox{iii) } \quad (a \wedge b) \wedge c = a \wedge (b \wedge c) \]
 for all $a,b,c$ in $E(n)$
and linear, i.e.
\[ \mbox{iv) } \quad a \wedge (\beta b + \gamma c) = \beta a \wedge 
b + \gamma a \wedge c \]
for all $\beta, \gamma$ in the field $\kr$ and $a,b,c$ in $E(n)$,
determine the algebra structure
on $E(n)$. For instance:
\begin{eqnarray*} 
e_5 \wedge (e_1 \wedge e_3) = &
e_5 \wedge e_1 \wedge e_3 & \\
= & -e_1 \wedge e_5 \wedge e_3  & \quad \mbox{(switch } e_5 \mbox{ and } e_1) \\
= &  e_1 \wedge e_3 \wedge e_5 & \quad  \mbox{(switch } e_5 \mbox{ and } e_3). \\
\end{eqnarray*}


\subsection{Abstract definition} Here we define the exterior
algebra using standard machinery from algebra.
Let $V$ be a vector space over $\kr$, and
denote by  $V^{\te p}$ the $p$-fold tensor product 
$V \te_{\kr} V \te_{\kr} \cdots \te_{\kr} V$. The free associative algebra on $V$
is the tensor algebra $T(V) = \oplus_{p\geq 0} V^{\te p}$ which comes
with the natural concatenation product
\[ (v_1 \te \cdots \te v_r) \cdot (w_1 \te \cdots \te w_s) = 
v_1 \te \cdots v_r \te w_1 \te \cdots \te w_s .\]
Let $R$ be the subspace of $V \te_{\kr} V$ generated by all
elements $v \te v$ where $v \in V$. The exterior algebra
is the quotient algebra of $T(V)$ by the relations $R$. 
More formally, let $\langle R \rangle$ be the two-sided ideal in $T(V)$
generated by  $R$. The exterior algebra $E(V)$ is the quotient algebra
$T(V)/\langle R \rangle$.
The product in this quotient algebra is commonly denoted by 
$\wedge$.  Let $e_1, \ldots, e_n$ be a basis for $V$.
We then have $e_i \wedge e_i = 0$ since $e_i \te e_i$ is a
relation in $R$. Similarly $(e_i + e_j)\wedge (e_i + e_j)$
is zero. Expanding this
\[ 0 =  e_i \wedge e_i + e_i \wedge e_j + e_j \wedge e_i + 
e_j \wedge e_j \]
we see that $e_i \wedge e_j = - e_j \wedge e_i$.
In fact we obtain $v \wedge w + w \wedge v = 0$ for any $v,w$ in $V$.
Hence when the characteristic of $\kr$ is not $2$, 
the exterior algebra may be defined
as $T(V) / \langle S_2 V \rangle$
where
\[ S_2 V = \{ v \te w + w \te v \,|\, v,w \in V\} \]
are the symmetric two-tensors in $V \te V$. 
The $p$'th graded piece of $E(V)$,
which is the image of $V^{\te p}$, is denoted as $\wedge^p V$.

\medskip
We shall in the following indicate:
\begin{itemize}
\item How central notions in various areas in mathematics arise
from natural structures on the exterior algebra.
\item How the exterior algebra or variations thereof is 
a natural tool in these areas.
\end{itemize}

\section{Combinatorics I:  Simplicial complexes and face rings}

\label{KombISec}
For simplicity denote the set $\{1,2, \ldots, n\}$ as $[n]$.
Each subset $\{i_1, \ldots, i_r\}$ of $[n]$ corresponds to a
monomial $e_{i_1} \wedge e_{i_2} \wedge \ldots \wedge e_{i_r}$
in the exterior algebra $E(n)$. 
For instance $\{2,5\} \sus [6]$ gives the monomial
$e_2 \wedge e_5$. It also gives the indicator vector
$(0,1,0,0,1,0) \in \hele_2^6$  (where $\hele_2 = \{0,1\}$),
with $1$'s at positions $2$ and $5$. We may then consider 
$e_2 \wedge e_5$ to have this multidegree.
This one-to-one correspondence between subsets of $[n]$ and
monomials in $E(n)$ 
suggests that it can be used
to encode systems of subsets of a finite set. The set systems 
naturally captured by virtue of $E(n)$ being an algebra, are
the {\it combinatorial simplicial complexes}. These are families
of subsets $\Delta$ of $[n]$ such that if $X$ is in $\Delta$
then any subset $Y$ of $X$ is also in $\Delta$.

\begin{example}  \label{KoIEksSicx}
Let $n = 6$. The sets 
\[ \{ 1,2 \}, \,\{3,4\}, \{3,5 \}, \, \{4,5,6 \} \]
together with all the subsets of each of these four sets
form a combinatorial simplicial complex.
\end{example}

The point of relating these to the algebra $E(n)$ is that combinatorial
simplicial complexes on $[n]$ are in one-to-one correspondence with 
$\hele_2^n$
graded ideals $I$ in $E(n)$ 
or equivalently with $\hele_2^n$-graded
quotient rings $E(n)/I$ of $E(n)$: To a simplicial complex $\Delta$
corresponds the monomial ideal $I_\Delta$ generated by
\[ \{ e_{i_1} \wedge \cdots \wedge e_{i_r} \,|\, \{i_1,\ldots, i_r \}
\not \in \Delta \}. \]
Note that the monomials $e_{i_1} \wedge \cdots \wedge e_{i_p}$ with
$\{i_1, \ldots, i_p\}$ in $\Delta$ then constitute a vector space basis for 
the quotient algebra $E(\Delta) = E(V)/I_\Delta$. We call this algebra the  
{\it exterior face ring}
of $\Delta$. 

For the simplicial complex in the example above, $E(\Delta)$ has
a basis:
\begin{itemize}
\item degree $0$: $1$,
\item degree $1$: $e_1, e_2, e_3, e_4, e_5, e_6$,
\item degree $2$: $e_1 \wedge e_2, e_3 \wedge e_4, e_3 \wedge e_5,
e_4 \wedge e_5, e_4 \wedge e_6, e_5 \wedge e_6$,
\item degree $3$: $e_4 \wedge e_5 \wedge e_6$. 
\end{itemize}

Although subsets $\{i_1, \ldots, i_r \}$ of $[n]$ most naturally
corresponds to monomials in $E(n)$, one can also consider the monomial
$x_{i_1} \cdots x_{i_r}$ in the 
polynomial ring $\kr[x_1, \ldots, x_n]$. (Note
however that monomials in this ring naturally corresponds to multisets rather
than sets.) If one to $\Delta$
associates the analog monomial ideal in this polynomial ring, the quotient
ring $\kr[\Delta]$ is the {\it Stanley-Reisner ring} or simply the 
{\it face ring} of $\Delta$. 

This opens up the arsenal of algebra to study $\Delta$. 
The study of $E(\Delta)$ and $\kr[\Delta]$ has particularly centered
around their minimal free resolutions and all the invariants
that arise from such. The study of $\kr[\Delta]$ was launched around 1975 
with a seminal paper by Hochster \cite{Ho} and Stanley's proof
of the Upper Bound Conjecture for simplicial spheres, see \cite{StCoCo}.
Although one
might say that $E(\Delta)$ is a more natural object associated to $\Delta$,
$\kr[\Delta]$ has been preferred for two reasons:
i) minimal free resolutions over $\kr[x_1, \ldots, x_n]$ are finite in contrast
to over the exterior algebra $E(n)$, 
ii) $\kr[\Delta]$ is commutative and the machinery 
for commutative rings is very well developed.

Since 1975 this has been a very active area of research with
various textbooks published \cite{StCoCo}, \cite{BHeCM},
\cite{MiStCoCo}, and \cite{HeHi}.
For the exterior face ring, see \cite{FlVaBi}.

\section{Topology}

Let $u_i = (0,\ldots,0,1,0, \ldots, 0)$ be the $i$'th unit
coordinate vector in $\Rn$.
To a subset $\{i_1, \ldots, i_r \}$ of $[n]$ 
we may associate the $(r-1)$-dimensional simplex which is
the convex hull of the points $u_{i_1}, \ldots, u_{i_r}$ in 
$\Rn$. For instance $\{2,3,5 \} \sus [6]$ gives the simplex
consisting of all points $(0,\la_2,\la_3,0,\la_5,0)$ where
$\la_i \geq 0$ and $\la_2 + \la_3 + \la_5 = 1$.

A combinatorial simplicial complex $\Delta$ has a natural topological
realization $X = |\Delta|$. It is the union of all the simplices
in $\Rn$ coming from sets $\{i_1, \cdots, i_r \}$ in $\Delta$. 

\begin{example} \label{TopEksSicx}
The simplicial complex given in Example \ref{KoIEksSicx}
has a topological realization which may be pictured as:

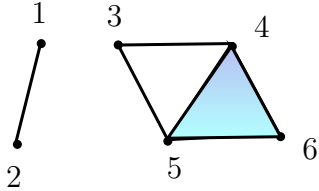
\scalebox{1} 
{
\begin{pspicture}(-2.5,-1.2729688)(2.9828124,1.2729688)
\definecolor{color56g}{rgb}{0.7254901960784313,0.7490196078431373,0.9490196078431372}
\definecolor{color56f}{rgb}{0.6941176470588235,1.0,1.0}
\psline[linewidth=0.04cm](0.6009375,0.71453124)(0.2609375,-0.6454688)
\psdots[dotsize=0.12](0.6009375,0.6745312)
\psdots[dotsize=0.12](0.2809375,-0.66546875)
\psline[linewidth=0.04](1.6009375,0.69453126)(2.2809374,-0.60546875)(3.1209376,0.65453124)(3.7609375,-0.52546877)(2.2809374,-0.60546875)
\psline[linewidth=0.04cm](1.6209375,0.65453124)(3.1609375,0.6745312)
\psdots[dotsize=0.12](1.6209375,0.65453124)
\psdots[dotsize=0.12](3.1409376,0.63453126)
\psdots[dotsize=0.12](2.2809374,-0.62546873)
\psdots[dotsize=0.12](3.7809374,-0.5654687)
\usefont{T1}{ptm}{m}{n}
\rput(0.57234377,1.0795312){$1$}
\usefont{T1}{ptm}{m}{n}
\rput(0.23234375,-1.1004688){$2$}
\usefont{T1}{ptm}{m}{n}
\rput(1.5723437,1.0395312){$3$}
\psline[linewidth=0.04,fillstyle=gradient,gradlines=2000,gradbegin=color56g,gradend=color56f,gradmidpoint=1.0](3.1409376,0.65453124)(2.2809374,-0.58546877)(3.7809374,-0.5654687)(3.1409376,0.6145313)(3.1609375,0.59453124)
\usefont{T1}{ptm}{m}{n}
\rput(3.5323439,0.89953125){$4$}
\usefont{T1}{ptm}{m}{n}
\rput(2.3723438,-0.98046875){$5$}
\usefont{T1}{ptm}{m}{n}
\rput(4.1723437,-0.72046876){$6$}
\end{pspicture} 
}

This is the disjoint union of a line segment and a disc
with a handle.
\end{example}

 We then say that
$\Delta$ gives a triangulation of the space $X$.
Now we equip $E(n)$ with a differential $d$ of degree $1$
by multiplying with $u = e_1 + e_2 + \cdots + e_n$.
Then $d(a) = u \wedge a$ and this 
is a differential since $d^2(a) = u \wedge u \wedge a = 0$.  
The monomial
$e_{i_1} \wedge e_{i_2} \wedge \cdots \wedge
e_{i_d}$ of degree $r$ is then mapped to
the degree $(d+1)$ sum:
\[ \sum_{i \not \in \{i_1, \ldots, i_r\}} e_i \wedge e_{i_1} \wedge \cdots
\wedge e_{i_d}. \]

The face ring $E(\Delta)$ is a quotient of $E(n)$ and
so we also obtain a differential on $E(\Delta)$. 
Letting $E(\Delta)^d$ be the degree $d$ part, this gives a complex
\[ E(\Delta)^0 \mto{d^0} E(\Delta)^1 \mto{d^1} E(\Delta)^2
\mto{d^2} \cdots . \]
From this complex and its dual we calculate the prime invariants in
topology, the cohomology and homology of the topological space $X$.
The cohomology is
\begin{equation*}
H^{i+1}(E(\Delta),d)  = \tilde{H}^i(X,\kr) \mbox{ for } i \geq 0,
\end{equation*}
where $\tilde{H}^i(X,\kr)$ is the reduced cohomology of $X$.
(For $i > 0$ this is simply the cohomology $H^i(X,\kr)$ while
for $i = 0$, this is the cokernel  $H^0(pt,\kr) \pil H^0(X,\kr)$. ) 
Dualizing the complex above we get $E(\Delta)^*$ as a subcomplex
of $E(n)^*$:
\[ \cdots \mto{\vardel_2} (E(\Delta)^*)_2 \mto{\vardel_1} 
(E(\Delta)^*)_1 \mto{\vardel_0} (E(\Delta)^*)_0. \]
Here $(E(\Delta)^*)_{d+1}$ has a basis consisting of monomials
\begin{equation} \label{TopLigFace} 
e^*_{i_0} \wedge \cdots \wedge e^*_{i_{d}} 
\end{equation}
where $\{i_0, \ldots, i_{d}\}$
are the $d$-dimensional faces of the simplicial complex $\Delta$.
The differential $\vardel$ is contraction with the
element $u$
\[ a \overset{\vardel}{\mapsto} u \neg a, \]
sending the monomial (\ref{TopLigFace}) to its boundary
\[ \sum_{j=0}^d (-1)^j e^*_{i_0} \wedge \cdots \wedge \widehat{e^*_{i_j}} \wedge
\cdots e^*_{i_d}. \]
(Here $\widehat{e^*_{i_j}}$ means omitting this term.)
The homology $H_i(E(\Delta)^*,\vardel)$  
computes the reduced simplicial homology of the space $X$.
In Example \ref{TopEksSicx} 
above we get $H_1(E(\Delta)^*, \vardel) = \kr$ one-dimensional,
one less than the number of components of $\Delta$, and 
$H_2(E(\Delta)^*, \vardel) = \kr$ 
one-dimensional since there is a non-contractible
$1$-cycle.

A good introduction to algebraic topology, starting from simplicial
complexes, is \cite{MuAT}.

\section{Lie theory} \label{LieSek}
Differential graded algebras (DGA) occur naturally in many areas. 
They provide the ``full story'' in contrast to graded
algebras which often are the cohomology of a DGA,
like the cochain complex of a topological space in contrast
to its cohomology ring. 

A DGA is graded algebra $A = \oplus_{p \geq 0} A_p$
with a differential $d$, i.e. $d^2 = 0$, which is a derivation,
i.e. for homogeneous elements $a,b$ in 
$A$ it satisfies:
\begin{equation} \label{LieLigDeriv}
 d(a \cdot b) = d(a)\cdot b + (-1)^{\deg(a)} a \cdot d(b), 
\end{equation}
The differential $d$ either has degree $1$ or $-1$ according
to whether it raises degrees by one or decreases degrees by one.

\medskip
What does it mean to give
a $\kr$-linear 
differential $d$ of degree $1$ on $E(V)$ such that $(E(V),d)$ becomes
a DGA?
Between degrees $1$ and $2$ we have a map 
\[ V \mto{d} \wedge^2 V.\]
By the definition of derivation \ref{LieLigDeriv}, it is easy to see
that any linear map between these vector spaces
extends uniquely to a derivation $d$ on
$E(V)$. 
Denote by $\frg$  the dual vector space $V^* = \Hom_k(V,k)$.
Dualizing the above map we get a map
\begin{eqnarray*} \wedge^2 \frg & \mto{d^*} &  \frg \\
x \wedge y & \mapsto &  [x,y].
\end{eqnarray*}
It turns out that $d$ gives a differential, i.e. $d^2 = 0$ if and only if
the map $d^*$ satisfies the Jacobi identity 
\[ [x,[y,z]] + [y,[z,x]] + [z,[x,y]] = 0. \]
Thus giving $E(V)$ the structure of a DGA with differential of degree $1$
is precisely equivalent to give $\frg = V^*$ the structure of a Lie algebra.

\medskip
The cohomology of the complex 
$(E(V),d)$ computes the {\it Lie algebra cohomology}
of $\frg$. If $\frg$ is the Lie algebra of a connected compact Lie
group $G$ (over $\kr = \mathbb{R}$), 
it is a theorem of Cartan 
see \cite[Cor.12.4]{Bre}, that
the cohomology ring $H^*(E(V),d)$ is isomorphic to the cohomology ring
$H^*(G, \RR)$.
A good and comprehensive introduction to Lie algebras is \cite{HuLie}.
The book \cite{FuHa} is much used as a reference book for the representations
of semi-simple Lie groups and Lie algebras. But there are so many
books on this, and the above are mentioned mostly because I learned from
these books.

Now denote the dual $V^*$ by $W$, and by $S(W) = \Symm(W)$ the symmetric
algebra, i.e. the polynomial ring whose variables are any basis of $W$.
The pair $E(V)$ and $S(W)$ is the prime example of a {\it Koszul dual} 
pair of algebras, see \cite{PolPos}, \cite{BGS} for the general framework
of Koszul duality.  
Furthermore, when we equip $E(V)$ with the differential $d$, 
the pair $(E(V),d)$ can be considered as the {Koszul dual}
of the enveloping algebra $U(\frg)$ of the Lie algebra $\frg$, 
see \cite{Pos}, \cite{FlKo} for the Koszul duality in the differential graded 
setting. Koszul duality gives functors
between the module categories of these algebras, 
which on suitable quotients of these
give an equivalence of categories.

\section{Combinatorics II: Hyperplane arrangements and 
the Orlik-Solomon algebra}
Simplicial complexes are basic combinatorial structures and
we have seen in Section \ref{KombISec} 
how they are captured by the exterior face ring
$E(\Delta)$. One
of the most successful unifying abstractions in combinatorics
is that of a {\it matroid} (a word giving more associations might
be {\it independence structures}) which is a special type of a simplicial 
complex.  To a matroid there is associated a quotient
algebra of the exterior algebra with remarkable connections to hyperplane
arrangements.

A prime source of matroids are from linear algebra.
Let us consider the vector space $\kr^m$ and let $x_1, \ldots, x_m$
be coordinate functions on this space. A linear form $v = \sum \la_j x_j$
gives a hyperplane in $\kr^m$:  the set of all points $(a_1, \ldots, a_m)$
in $\kr^m$ such that $\sum_j \la_j a_j = 0$. 
A set of linear forms $v_1, \ldots, v_n$
determines hyperplanes $H_1, H_2, \ldots, H_n$. We call this a 
{\it hyperplane arrangement}. 
It turns out that a number of essential properties of the
hyperplane arrangement is determined by the linear dependencies of the
linear forms $v_1, \ldots, v_n$. 
We get a combinatorial simplicial complex $M$ on $[n]$ consisting of
all subsets $\{i_1, \ldots, i_r\}$ of $[n]$ such that $v_{i_1}, \ldots, v_{i_r}$
are linearly independent vectors. But there is more structure on this $M$,
making it a matroid. More precisely a simplicial complex $M$
on $[n]$ is a {\it matroid} if the following extra condition holds.

\medskip
{\it If $X$ and $Y$ are independent sets of $M$, 
with the cardinality of $Y$ larger than that of $X$,
there is $y \in Y \backslash X$ such that $X \cup \{y \}$ is independent.
}

\medskip
The elements of $M$ are called the {\it independent sets} of the
matroid, while subsets of $[n]$ not in $M$ are {\it dependent}. 
The diversity which the abstract notion of a matroid captures
is illustrated by the following examples where we give independent
sets of matroids:
\begin{itemize}
\item Linear independent subsets of a set of vectors 
$\{v_1,v_2, \ldots, v_n\}$.
\item Edge sets of graphs which do not contain a cycle.
\item Partial transversals of a family of sets $A_1, A_2, \ldots, A_N$.
\end{itemize}

\begin{example}
Consider the hyperplane arrangement in $\CC^2$ given by
the two coordinate functions $v_1 = x_1$ and $v_2 = x_2$. 
The complement $\CC^2 \backslash H_1 \cup H_2$ consists of the
pairs $(a,b)$ with non-zero coordinates, i.e. the complement
is $(\CC^*)^2$ where $\CC^* = \CC \backslash \{0 \}$. Since
$\CC^*$ is homotopy equivalent to the circle $S^1$, 
the complement $\CC^2 \backslash H_1 \cup H_2$ is 
homotopy equivalent to the torus $S^1 \times S^1$. The cohomology
ring of this torus is the exterior algebra $E(2)$. 
\end{example}

This example generalizes to a description of the cohomology 
ring of the complement of any hyperplane arrangement in $\CC^m$. 
The matroid $M$ of the hyperplane arrangment, being a simplicial 
complex, gives by Section \ref{KombISec} a monomial ideal $I_M$ in $E(n)$.
The dual element $u = e_1^* + \cdots + e_n^*$ gives a contraction
$a \mto{\vardel} u \neg a$ sending $e_{i_1} \wedge \cdots \wedge e_{i_r}$
to
\[ \sum_j (-1)^ j e_{i_1} \wedge \cdots \hat{e_{i_j}} \cdots \wedge e_{i_r}
\]
(here $\hat{e_{i_j}}$ means omitting this term). 
Now $I_M + 
\vardel(I_M)$ also becomes an ideal in $E(n)$. The quotient
$A(M) = E(n)/(I_M + \vardel(I_M))$ is the
{\it  Orlik-Solomon algebra} associated to the hyperplane arrangement.
In 1980 Peter Orlik and Louis Solomon proved the following amazing result,
\cite{OrSo}.

\begin{theorem}
Let $T = \CC^m \backslash \cup_i H_i$ be the complement of a complex
hyperplane arrangement. Then the algebra $A(M)$ is the 
cohomology ring $H^*(T, \kr)$. 
\end{theorem}

\begin{example}
Consider $u_1 = x_1 - x_2$, $u_2 = x_2 - x_3$ and $u_3 = x_3 - x_1$
on $\CC^3$. There is one dependency her, between $u_1, u_2$ and $u_3$.
Thus the Orlik-Solomon algebra is $E(3)$ divided by the relation 
\[ \vardel (e_1 \wedge e_2 \wedge e_3) = e_1 \wedge e_2 - e_1 \wedge e_3
+ e_2 \wedge e_3. \]
The quotient algebra of $E(3)$ by these relations has dimensions $1,3$ and $2$ 
in degrees $0, 1$ and 
$2$ respectively. For the complement $T = \CC^3 \backslash H_1 \cup H_2 \cup H_3$
we therefore have 
\[ H^0(T, \CC) = \CC, \quad H^1(T, \CC) = \CC^3, H^2(T, \CC) = \CC^2. \]
\end{example}

For hyperplane arrangements or more generally for matroids, the 
Orlik-Solomon algebra has been much studied during the later years,
\cite{YuzOS}, \cite{FaOS}. The algebraic
properties of the Orlik-Solomon algebra gives a number of natural 
invariants for hyperplane arrangements.

\section{Mathematical physics}

The Clifford algebra may be viewed as a deformation of the exterior 
algebra. The exterior algebra $E(V)$ is defined as the
quotient algebra $T(V)/\langle S_2 V \rangle$.
Fix now a symmetric bilinear form $b : S_2 V \pil k$.
Let $R = \{ r - b(r) \, | \, r \in S_2 V \}$. 
The Clifford algebra
is the quotient of the tensor algebra by the relations $R$: 
\[ Cl_b = T(V)/ \langle R \rangle . \]
Like the exterior algebra $E(V)=E(n)$ it has a basis consisting of all products
$e_{i_1}\cdot \cdots \cdot e_{i_r}$ for subsets $\{i_1 < \cdots < i_r\} $
of $\{1, \ldots, n\}$, and so is also of dimension
$2^n$. Note that we get the exterior
algebra when the bilinear form $b = 0$. 

Clifford algebras are mostly applied when 
the field $\kr$ is the real numbers $\RR$. 
Let $V = \langle i \rangle $ be a
one dimensional vector space, and the quadratic form be given by
$i^2 \mtoelm{b} -1$. Then the associated Clifford algebra is the complex
numbers. When $V = \langle i, j \rangle$ is a two-dimensional space and
\[ i \te i  \mtoelm{b} -1,\quad i \te j + j \te i \mtoelm{b} 0, \quad 
j \te j \mtoelm{b} -1, \]
we obtain the quaternions. 
In general for a real symmetric form $b$, we may find a basis 
for  $V$ such that if 
$x_1, \ldots, x_n$ are the coordinate functions, the form is
\[ \sum_{i=1}^p x_i^2 - \sum_{i= p+1}^{p+q} x_i^2.\]
Such a Clifford algebra is denoted $Cl_{p,q}$. So $Cl_{0,1}$ is
the complex numbers and $Cl_{0,2}$ is the quaternions.
Clifford algebras have interesting periodic behaviour:
$Cl_{p+1,q+1}$ is isomorphic to the $2 \times 2$-matrices
$M_2(Cl_{p,q})$, and each of 
$Cl_{p+8,q}$ and $Cl_{p,q+8}$ is isomorphic to the $16 \times 16$-matrices
$M_{16}(Cl_{p,q})$. Thus Clifford algebras over the reals are
essentially classified by $Cl_{p,0}$ and $Cl_{0,p}$ for $p \leq 7$.
For an introduction to Clifford algebras, see \cite{GaInCl}.

When $Cl_{p,q}$ is a simple algebra, and $Cl_{p,q} \pil End(W)$
is an irreducible representation of $Cl_{p,q}$, then $W$ is called
a spinor space. 
These representations occur a lot
in mathematical physics. For instance $Cl_{3,1}$ is isomorphic
to $M_4(\RR)$ and this representation on $\RR^4$ is the Minkowski
space with one time-dimension and three space dimensions.
A pioneer in the application of Clifford algebras in mathematical
physics is David Hestenes, \cite{HeSpace}, \cite{HeGeAl}, \cite{HeNeMe},
where he envisions the complete use of it in classical mechanics. He calls
this geometric algebra. 
The book \cite{DoGeAlFy} offers a leisurely introduction to 
the application of geometric algebra in physics.
Basil Hiley is another 
advocate for the algebraic approach to quantum mechanics, \cite{HiCl}:

\begin{quote}
that quantum phenomena per se can be entirely described in terms of 
Clifford algebras taken over the reals without the need to appeal to specific 
representations in terms of wave functions in a Hilbert space. This removes 
the necessity of using Hilbert space and all the physical imagery that goes 
with the use of the wave function.
 \end{quote}

\section{Algebraic geometry}
Finitely generated graded modules over exterior algebras seem
far removed from geometry. However we shall see that they encode
perhaps the most significant invariants of algebraic geometry, the cohomological
dimensions of twists of sheaves on projective spaces.

\begin{example} \label{AlgEksTate}
Let $n = 2$ and $E = E(2)$.
Consider the map of free $E$-modules
\[ E \mto{ d = \left [ \begin{matrix} e_2 \\ e_1 \end{matrix} \right ] } E^2.
\]
Writing $E^2 = Eu_1 \oplus Eu_2$, 
the cokernel of this map is a module 
$M = Eu_1 \oplus Eu_2 / \langle e_2 u_1 + e_1 u_2 \rangle$. 
Such a map may, as we shall shortly explain, be
completed to a complex of free $E$-modules (we let $d^0 = d$):
\begin{eqnarray} \label{AlgLigTateP1}
\cdots \pil E^2  \xrightarrow{d^{-2} = \left [ \begin{matrix} e_1 & e_2 \end{matrix}
\right ]}
& E \xrightarrow{ d^{-1} = \left [ \begin{matrix} e_1 \wedge e_2 \end{matrix}  \right ]} E   \xrightarrow{ d^0 = \left [ \begin{matrix} e_2 \\ e_1 \end{matrix} \right ] } \\ \notag
& E^2 
\xrightarrow{ d^1 = \left [ \begin{matrix} e_2 & 0 \\ e_1 & e_2 \\ 0 & e_1 \end{matrix} \right ] }
E^3 \pil E^4 \pil \cdots 
\end{eqnarray}

It is a complex since $d^p \circ d^{p-1} = 0$, as is easily verified. 
It is also exact at each place, meaning that the kernel of $d^p$ 
equals the image of $d^{p-1}$ for each $p$. Thus it is an {\it acyclic
complex}.
\end{example}

Every finitely generated graded module $M$ over $E(V)$,
or equivalently map $d$, gives rise to such an acyclic complex. In the
example above the ranks of the free modules follow a simple pattern
$1,2,3,4,\cdots$. But in general, what are these
ranks? Can they be given a meaningful interpretation? Indeed, a
fairly recent discovery from 2003, \cite{EFS}, tells us this is the case.

\medskip
 As a module over itself $E(V)$ is both a projective 
and an injective module. Given a finitely generated graded module $M$ over
$E(V)$ one can make a minimal free (and so projective) resolution
\[ P^\dt \pil M, \quad \mbox{where } 
P^p = \oplus_{q \in \hele} W^p_q \te_k E, \]
(the $W^p_q$ are vector spaces whose elements have degree $q$)
and a minimal injective resolution
\[ M \pil I^{\dt}, \quad \mbox{where }
I^p = \oplus_{q \in \hele} W_q^p \te_k E \]
and splice these together to an acyclic complex (like in Example 
\ref{AlgEksTate} above)
\begin{equation} \label{AlgeoLigTateres} 
T : \cdots \pil P^{-1} \pil P^0 \pil I^1 \pil I^2 \pil \cdots 
\end{equation}
called the {\it Tate resolution} of $M$. So we get a correspondence:
\begin{equation} \label{AlgeoLigTate} \text{graded modules over } E(V) \leadsto
\text{Tate resolutions over } E(V).
\end{equation}

\medskip
Now let us pass to another construction starting from the
finitely generated graded module $M = \oplus_{i=a}^b M_i$
over the exterior algebra $E(V)$. Let $W$ be the dual space $V^*$. 
The multiplication
$V \te_k M_i \pil M_{i+1}$ gives a map
\begin{equation} 
M_i \pil W \te_k M_{i+1} \label{AlgLigMW}
\end{equation}
Let $S = \Symm(W)$ be the symmetric algebra.
The maps (\ref{AlgLigMW}) gives rise to maps
\begin{equation} \label{AlggeoLigLin} \cdots \pil S \te_k M_p \mto{d^p} S \te_k M_{p+1} \mto{d^{p+1}} \cdots 
\mto{d^{q-1}} S \te_k M_q \pil \cdots.
\end{equation}
These maps give a (bounded) complex of $S$-modules, 
i.e. $d^{i+1} \circ d^i = 0$. (A correspondence within the framework
of Koszul duality, 
mentioned in Section \ref{LieSek}.) Any finitely generated graded $S = \Symm(W)$-module may be {\it sheafified} to a
{\it coherent sheaf} on the projective space $\PW$. 
In particular we may sheafify the above complex and get a complex
of coherent sheaves:
\begin{align} \label{AlgeoLigBGG} & \text{graded modules over } E(V) \\ \notag \leadsto 
& \text{ bounded complexes of coherent sheaves on } \PW. 
\end{align}
This correspondence is
from 1978, \cite{BGG}, 
and is the celebrated Bernstein-Gelfand-Gelfand (BGG) correspondence. 
Somewhat more refined it may be described as an equivalence of 
categories between suitable categories of the objects in (\ref{AlgeoLigBGG}).

\medskip
The amazing thing is that if $M$ via the BGG-correspondence 
(\ref{AlgeoLigBGG}) 
gives a coherent sheaf $\gF$ on 
$\PW$ (this means that the sheafification of the complex (\ref{AlggeoLigLin})
has only one nonzero cohomology sheaf $\gF$), then we can read off all the 
{\it sheaf cohomology}
groups of all twists of $\gF$ from the Tate resolution $T$ which we get 
via the correspondence
(\ref{AlgeoLigTate}).

\begin{theorem}[{\cite[Thm.4.1]{EFS}}] 
If $M$ via the BGG-correspondence (\ref{AlgeoLigBGG}) 
gives a coherent sheaf $\gF$ on the projective space $\PW$, and the 
Tate resolution of $M$ is (\ref{AlgeoLigTateres}), then
the sheaf cohomology
\[ H^p (\PW, \gF(q)) = W_q^{p+q}. \]
\end{theorem}

Returning to the initial example in this section, the sheaf
corresponding to this module $M$ is the structure
sheaf $\gO_{\PP^1}$ of the projective line $\PP^1 = \PP(W)$. 
The Tate resolution (\ref{AlgLigTateP1}) 
therefore tells us that:
\[ H^0(\PP^1, \gO_{\PP^1}(d)) = \begin{cases} \kr^{d+1}, &  d \geq 0 \\
 0, & d < 0 \end{cases}, \quad
H^1(\PP^1, \gO_{\PP^1}(-d)) = \begin{cases} \kr^{d-1}, & -d < 0 \\
0, & -d \geq 0 \end{cases}.
\]

\medskip
An important feature of a Tate resolution 
$T$ is that it is fully determined by an
arbitrary differential $T^i \mto{d^i} T^{i+1}$. This is because
$T^{\leq i}$ is a minimal projective
resolution of $\im d^i$ and $T^{> i}$ is a minimal injective resolution
of $\im d^i$. 
This gives us incredible freedom in construction.
An arbitrary homogeneous matrix $A$
of exterior forms gives a map 
\begin{equation} \oplus_q W_q^0 \te_k E \mto{d_A} \oplus_q W_q^1 \te_k E. 
\label{AlggeoLigA}
\end{equation}
The module $M = \im d_A$ over the exterior algebra 
then gives a complex of coherent sheaves $\gF^\dt$
on $\PW$ by (\ref{AlgeoLigBGG}), 
and all such bounded complexes on $\PW$ do (in a suitable sense) 
come from such a homogeneous matrix
$A$ of exterior forms. 
Thus bounded complexes
of coherent sheaves on projective spaces can be specified by giving 
a homogeneous matrix $A$ 
of exterior forms, and any matrix $A$ will give some such bounded complex.
The Tate resolution associated to $M$ and  $\gF^\dt$ is the complex we 
obtain by taking
a minimal projective resolution of $\ker d_A$ in (\ref{AlggeoLigA}) and
a minimal injective resolution of $\coker d_A$, and this
resolution tells us the cohomology of the complex of coherent sheaves.

\begin{example} Let $V$ the five-dimensional vector space generated
by $\{ e_1, e_2, e_3, e_4, e_5\}$. The matrix
\[ A = \begin{bmatrix} e_1 \wedge e_2 & e_2 \wedge e_3 & e_3 \wedge e_4 &
e_4 \wedge e_5 & e_5 \wedge e_1 \\
e_3 \wedge e_5 & e_4 \wedge e_1 & e_5 \wedge e_2 & e_1 \wedge e_3 &
e_2 \wedge e_4 \end{bmatrix} \]
gives a map $E^5 \pil E^2$.
Via the BGG-correspondence (\ref{AlgeoLigBGG})
this gives the celebrated Horrocks-Mumford
bundle on $\mathbb{P}^4$ discovered forty years ago, 
\cite[Section 8]{EFS}. 
In characteristic zero this is essentially the 
only known indecomposable rank two bundle on any projective space 
of dimension greater or equal to four. It is an intriguing problem
to use the methods above to try to construct new bundles of rank
$\leq n-2$ on a projective space $\mathbb{P}^n$, but to our knowledge
nobody has yet been successful.
\end{example}

Tate resolutions and algebraic geometry are treated in the books
\cite{EiGeSy} and \cite{M2Co}. The software program \cite{M2} contains
the package BGG for doing computations with Tate resolutions.

\section{Modelling and computations}
The last fifteen years have seen a flurry of books and
treatises giving applications of exterior algebras and Clifford algebras,
usually under the name "geometric algebra".
Groups at the University of Cambridge and the University of Amsterdam
have been particularly active in promoting geometric algebra.
The book "Geometric algebra for physicists" \cite{DoGeAlFy}
by C.Doran and A.Lasenby, 2003 is a very well written and readable
introduction to exterior algebras, Clifford algebras, and their
applications in all areas of physics, following the ideas outlined
by D.Hestenes. A more advanced 
treatment is \cite{BaCl}.
The book "Geometric algebra for computer scientists: An object oriented
approach", \cite{DoGeAlCo} by L.Dorst, D.Fontijne, and S.Mann
shows geometric algebra as an effective tool to describe a
variety of geometric models, involving linear spaces, circles, spheres, 
rotations,
and reflections. In particular it considers the conformal geometric
model developed in \cite{HeConform}. 
"Geometric algebra for engineers", \cite{Per} by C.Perwass similarly
applies geometric algebra to models occurring in engineering: camera
positions, motion tracking, and statistics. It 
also considers numerical aspects of its implementation.
Other books on geometric algebra and its use in computer modelling
and engineering are
\cite{HeOld}, \cite{HiFo}, \cite{SoGeAlg}, 
\cite{DoDoGeAl}, \cite{ViGeAlCo}, and \cite{ViCoAn}.
The book \cite{DoLaGeAlPrac} gives a panorama of applications by a
wide range of authors.

The comprehensive book "Grassmann algebra" 
\cite{Browne} considers all aspects of computations concerning the 
exterior algebra with Mathematica. It treats the exterior, interior,
and regressive products and geometric interpretations.
A second volume treats the generalized Grassmann product
which constitute an intermediate chain of products between the
exterior and interior products, and applications to hypercomplex numbers 
and to mechanics.
Other treatises with a more purely mathematical focus are
\cite{Schulz} and \cite{LuSv}.

\bibliographystyle{amsplain}
\bibliography{Bibliography}

\providecommand{\bysame}{\leavevmode\hbox to3em{\hrulefill}\thinspace}
\providecommand{\MR}{\relax\ifhmode\unskip\space\fi MR }
\providecommand{\MRhref}[2]{%
  \href{http://www.ams.org/mathscinet-getitem?mr=#1}{#2}
}
\providecommand{\href}[2]{#2}
\begin{thebibliography}{10}

\bibitem{BaCl}
William~E Baylis, \emph{Clifford (geometric) algebras: With applications in
  physics, mathematics, and engineering}, Springer, 1996.

\bibitem{BGS}
Alexander Beilinson, Victor Ginzburg, and Wolfgang Soergel, \emph{Koszul
  duality patterns in representation theory}, Journal of the American
  Mathematical Society \textbf{9} (1996), no.~2, 473--527.

\bibitem{BGG}
I.~N. Bern{\v{s}}te{\u\i}n, I.~M. Gel'fand, and S.~I. Gel'fand,
  \emph{{Algebraic vector bundles on {${\bf P}^{n}$} and problems of linear
  algebra}}, Funktsional. Anal. i Prilozhen. \textbf{12} (1978), no.~3, 66--67.

\bibitem{Bre}
Glen~E Bredon, \emph{Topology and geometry}, vol. 139, Springer Science \&
  Business Media, 1993.

\bibitem{Browne}
John Browne, \emph{Grassman algebra}, vol.~1, Create Space Independent
  Publishing Platform, 2012.

\bibitem{BHeCM}
Winfried Bruns and J{\"u}rgen Herzog, \emph{Cohen-macaulay rings}, Cambridge
  University Press, 1998.

\bibitem{Crowe}
Michael~J Crowe, \emph{A history of vector analysis: The evolution of the idea
  of a vectorial system}, Courier Dover Publications, 1967.

\bibitem{DoGeAlFy}
Christian Doran and Anthony Lasenby, \emph{Geometric algebra for physicists},
  Cambridge University Press, 2007.

\bibitem{DoDoGeAl}
Leo Dorst, Chris Doran, and Joan Lasenby, \emph{Applications of geometric
  algebra in computer science and engineering}, Springer, 2002.

\bibitem{DoGeAlCo}
Leo Dorst, Daniel Fontijne, and Stephen Mann, \emph{Geometric algebra for
  computer science (revised edition): An object-oriented approach to geometry},
  Morgan Kaufmann, 2009.

\bibitem{DoLaGeAlPrac}
Leo Dorst and Joan Lasenby, \emph{Guide to geometric algebra in practice},
  Springer, 2011.

\bibitem{M2Co}
D.~Eisenbud and et. al., \emph{{Computations in algebraic geometry with
  Macaulay 2}}, vol.~8, Springer, 2002.

\bibitem{EFS}
D.~Eisenbud, G.~Fl{\o}ystad, and F.O. Schreyer, \emph{{Sheaf cohomology and
  free resolutions over exterior algebras}}, Transactions of the American
  Mathematical Society \textbf{355} (2003), no.~11, 4397--4426.

\bibitem{EiGeSy}
David Eisenbud, \emph{{The Geometry of Syzygies: A Second Course in Algebraic
  Geometry and Commutative Algebra}}, vol. 229, Springer, 2005.

\bibitem{FaOS}
Michael Falk, \emph{{Combinatorial and algebraic structure in Orlik--Solomon
  algebras}}, European Journal of Combinatorics \textbf{22} (2001), no.~5,
  687--698.

\bibitem{FlVaBi}
G.~Fl{\o}ystad and J.E. Vatne, \emph{{(Bi-) Cohen--Macaulay simplicial
  complexes and their associated coherent sheaves}}, Communications in algebra
  \textbf{33} (2005), no.~9, 3121--3136.

\bibitem{FlKo}
Gunnar Fl{\o}ystad, \emph{Koszul duality and equivalences of categories},
  Transactions of the American Mathematical Society \textbf{358} (2006), no.~6,
  2373--2398.

\bibitem{FuHa}
W.~Fulton and J.~Harris, \emph{{Representation theory: A first course}}, GTM
  129, Springer, 1991.

\bibitem{GaInCl}
D.J.H. Garling, \emph{{Clifford algebras: An introduction}}, vol.~78, Cambridge
  University Press, 2011.

\bibitem{Gra}
Hermann Grassmann and Lloyd~C Kannenberg, \emph{Extension theory}, American
  Mathematical Society Providence, 2000.

\bibitem{M2}
Daniel~R. Grayson and Michael~E. Stillman, \emph{Macaulay2, a software system
  for research in algebraic geometry}, Available at
  http://www.math.uiuc.edu/Macaulay2/.

\bibitem{HeHi}
J{\"u}rgen Herzog and Takayuki Hibi, \emph{Monomial ideals}, Springer, 2011.

\bibitem{HeSpace}
David Hestenes, \emph{Space-time algebra}, Gordon and Breach, 1966.

\bibitem{HeGeAl}
\bysame, \emph{{Clifford algebra to geometric algebra: A unified language for
  mathematics and physics}}, D.Reidel publishing company, 1984.

\bibitem{HeNeMe}
\bysame, \emph{New foundations for classical mechanics}, Springer, 1999.

\bibitem{HeOld}
\bysame, \emph{Old wine in new bottles: A new algebraic framework for
  computational geometry}, Geometric Algebra with Applications in Science and
  Engineering, Springer, 2001, pp.~3--17.

\bibitem{HiFo}
Dietmar Hildenbrand, \emph{Foundations of geometric algebra computing}, vol.~8,
  Springer, 2012.

\bibitem{HiCl}
B.J. Hiley and R.E. Callaghan, \emph{{Clifford algebras and the Dirac-Bohm
  quantum Hamilton-Jacobi equation}}, Foundations of Physics \textbf{42}
  (2012), no.~1, 192--208.

\bibitem{Ho}
Melvin Hochster, \emph{Cohen-{M}acaulay rings, combinatorics, and simplicial
  complexes}, Ring theory, {II} ({P}roc. {S}econd {C}onf., {U}niv. {O}klahoma,
  {N}orman, {O}kla., 1975), Dekker, New York, 1977, pp.~171--223. Lecture Notes
  in Pure and Appl. Math., Vol. 26.

\bibitem{HuLie}
James~E Humphreys, \emph{{Introduction to Lie algebras and representation
  theory}}, vol. 1980, Springer New York, 1972.

\bibitem{HeConform}
Hongbo Li, David Hestenes, and Alyn Rockwood, \emph{Spherical conformal
  geometry with geometric algebra}, Geometric computing with Clifford algebras,
  Springer, 2001, pp.~61--75.

\bibitem{LuSv}
Douglas Lundholm and Lars Svensson, \emph{Clifford algebra, geometric algebra,
  and applications}, 2009.

\bibitem{MiStCoCo}
Ezra Miller and Bernd Sturmfels, \emph{Combinatorial commutative algebra}, vol.
  227, Springer, 2004.

\bibitem{MuAT}
James~R Munkres, \emph{Elements of algebraic topology}, vol.~2, Addison-Wesley
  Reading, 1984.

\bibitem{OrSo}
Peter Orlik and Louis Solomon, \emph{Combinatorics and topology of complements
  of hyperplanes}, Invent. Math. \textbf{56} (1980), no.~2, 167--189.

\bibitem{Per}
Christian Perwass, \emph{Geometric algebra with applications in engineering},
  vol.~4, Springer, 2008.

\bibitem{Petsche}
Hans-Joachim Petsche, Albert~C Lewis, J{\"o}rg Liesen, and Steve Russ (eds.),
  \emph{{From past to future: Gra{\ss}mann's work in context: Gra{\ss}mann
  Bicentennial Conference, September 2009}}, Springer, 2010.

\bibitem{PolPos}
Alexander Polishchuk and Leonid Positselski, \emph{Quadratic algebras},
  University Lecture Series, vol.~37, 2005.

\bibitem{Pos}
Leonid~Efimovich Positsel'skii, \emph{Nonhomogeneous quadratic duality and
  curvature}, Functional analysis and its applications \textbf{27} (1993),
  no.~3, 197--204.

\bibitem{Schub}
Gert Schubring (ed.), \emph{{Hermann G\"unther Gra{\ss}mann (1809-1877):
  visionary mathematician, scientist and neohumanist scholar}}, Kluwer Academic
  Publishers, 1997.

\bibitem{Schulz}
William~C Schulz, \emph{Theory and application of grassmann algebra}, Available
  at http://www.cefns.nau.edu/~schulz/grassmann.pdf, 2011, Preprint.

\bibitem{SoGeAlg}
Gerald Sommer, \emph{Geometric computing with clifford algebras: theoretical
  foundations and applications in computer vision and robotics}, Springer,
  2001.

\bibitem{StCoCo}
Richard~P Stanley, \emph{Combinatorics and commutative algebra}, Birkh{\"a}user
  Boston, 2004.

\bibitem{ViGeAlCo}
John~A Vince, \emph{Geometric algebra for computer graphics}, vol.~1, Springer,
  2008.

\bibitem{ViCoAn}
\bysame, \emph{Geometric algebra: An algebraic system for computer games and
  animation}, Springer, 2009.

\bibitem{YuzOS}
S.A. Yuzvinsky, \emph{{Orlik--Solomon algebras in algebra and topology}},
  Russian Mathematical Surveys \textbf{56} (2007), no.~2, 293.

\end{thebibliography}

\end{document}